\documentclass[11pt]{article}
\usepackage[english]{babel}
\usepackage{amssymb}
\usepackage{amsmath}
\usepackage{graphicx}
\usepackage{graphics}
\usepackage{amsthm}

\input pictex.tex

\newcommand{\clo}{\mathrm{S}^1}

\newcommand{\esp}{\thinspace}
\theoremstyle{definition}

\newtheorem{thm}{Theorem}[section]

\newtheorem{lem}[thm]{Lemma}

\input epsf

\begin{document}

\date{}
\author{Victor Kleptsyn \hspace{0.7cm} \& \hspace{0.7cm} Andr\'es Navas}

\title{A Denjoy Theorem for commuting circle 
diffeomorphisms with mixed H\"older derivatives}
\maketitle

\vspace{-0.25cm}

\noindent{\bf Abstract.} We prove that if $d \geq 2$ is an integer number 
and $f_k, \esp k \in \{1,\ldots,d\},$ are $C^{1+\tau_k}$ commuting circle 
diffeomorphisms with $\tau_k \!\! \in ]0,1[$ and $\tau_1 + \cdots + \tau_d > 1$, 
then the $f_k$'s are simultaneously (topologically) conjugate to rotations 
provided that their rotation numbers are independent over the rationals.

\vspace{0.15cm}

\noindent{\bf Keywords:} Denjoy Theorem, centralizers, 
H\"older derivative.


\markright{\sc Denjoy Theorem for commuting diffeomorphisms}

\vspace{0.3cm}

\section*{Introduction}

\hspace{0.5cm} Starting from the seminal works by Poincar\'e \cite{poincare} and 
Denjoy \cite{denjoy}, a deep theory for the dynamics of circle diffeomorphisms 
has been developed by many authors \cite{arnold,herman,KO,yoccoz}, and most 
of the fundamental related problems have been already solved. Quite 
surprisingly, the case of several commuting diffeomorphisms is rater 
special, as it was pointed out for the first time by Moser \cite{moser} 
in relation to the problem of the smoothness for the simultaneous 
conjugacy to rotations. Roughly speaking, in this 
case it should be enough to assume a joint Diophantine condition on the rotation 
numbers which does not imply a Diophantine condition for any of them (see the recent 
work \cite{fayad} for the solution of the $C^{\infty}$ case of Moser's problem).

A similar phenomenon concerns the classical Denjoy Theorem. Indeed, in \cite{DKN} 
it was proved that if $d \geq 2$ is an integer number and $\tau > 1/d$, then the 
elements $f_1,\ldots,f_d$ of any family of $C^{1+\tau}$ commuting circle 
diffeomorphisms are simultaneously (topologically) conjugate to rotations 
provided that their rotation numbers are independent over the rationals 
(that is, no non trivial linear combination of them with rational 
coefficients equals a rational number). In other words, the classical (and 
nearly optimal) $C^2$ hypothesis for Denjoy Theorem can be weakened in the 
case of several commuting diffeomorphisms. The first and main result of this 
work is a generalization of this fact to the case of different regularities.

\vspace{0.4cm}

\noindent{\bf Theorem A.} {\em Let $d \geq 2$ be an integer number 
and $\tau_1,\ldots,\tau_d$ be real numbers in $]0,1[$ such that 
$\tau_1 + \cdots + \tau_d > 1$. If $f_k, \esp k \in \{1,\ldots,d\},$ 
are respectively $C^{1+\tau_k}$ circle diffeomorphisms which have 
rotation numbers independent over the rationals and which do commute, 
then they are simultaneously (topologically) conjugate to rotations.}
\vspace{0.4cm}

Since the probabilistic arguments of \cite{DKN} cannot be  
applied to the case of different regularities, the preceding 
result is much more than a straightforward generalization of 
Theorem A of \cite{DKN}. Indeed, for the proof here we use 
a key new argument which is somehow more deterministic. 

Theorem A is (almost) optimal (in the H\"older scale), in the sense that if one 
decreases slightly the regularity assumptions then it is no longer true. The 
following result relies on classical constructions by Bohl \cite{bol}, Denjoy 
\cite{denjoy}, Herman \cite{herman}, and Pixton \cite{Pix}, and its proof 
consists on an easy extension of the construction given by Tsuboi in \cite{TsP}.

\vspace{0.4cm}

\noindent{\bf Theorem B.} {\em Let $d \geq 2$ be an integer number 
and $\tau_1,\ldots,\tau_d$ be real numbers in $]0,1[$ such that 
$\tau_1 + \cdots + \tau_d < 1$. If $\rho_1,\ldots,\rho_d$ are 
elements in $\mathbb{R}/\mathbb{Z}$ which are independent 
over the rationals, then there exist $C^{1+\tau_k}$ circle 
diffeomorphisms $f_k, \esp k \in \{1,\ldots,d\}$, having rotation 
numbers $\rho_k$, which do commute, and such that none of them is 
topologically conjugate to a rotation.}

\vspace{0.4cm}

It is well known that the techniques developed for Denjoy Theory can be applied 
to the study of group actions on the interval. In this direction we should point out 
that the methods of this paper also allow to extend (in a straightforward way) the 
so called ``Generalized Kopell Lemma" and the 
``Denjoy-Szekeres Type Theorem" (Theorems B and C of 
\cite{DKN} respectively) for Abelian groups of interval diffeomorphisms under 
analogous hypothesis of different regularities. Furthermore, the construction 
of counter-examples for both of them when these hypothesis do not hold can be 
also extended to this context. We leave the verification of all of this to the reader.


\vspace{0.35cm}

\noindent{\bf Acknowledgments.} It is a pleasure to thank Bassam Fayad 
and Sergey Voronin for their encouragements, as well as the Independent 
University of Moscow for the hospitality during the conference ``Laminations 
and Group Actions in Dynamics" held in February 2007. The first author was 
supported by the Swiss National Science Foundation. This work was also funded 
by the RFBR grants 7-01-00017-a and CNRS-L$_{-}$a 05-01-02801, and by the 
CONICYT grant 7060237. 


\section{A general principle revisited}

\hspace{0.5cm} As it is well known since the classical works by Denjoy, Schwartz 
and Sacksteder \cite{denjoy,Sa,Sc}, if $I$ is a wandering interval\footnote{We say 
that an interval is {\em wandering} if its images by different elements of the 
underlying semigroup are disjoint.} for the dynamics of a finitely generated 
semigroup $\Gamma$ of $C^{1+lip}$ diffeomorphisms of the 
closed interval or the circle (on which we will 
always consider the normalized length), one can control 
the distortion of the elements of $\Gamma$ over (a slightly larger 
interval than) $I$ in terms of the sum of the lengths 
of the images of $I$ along the corresponding sequence 
of compositions and a uniform Lipschitz constant for the 
derivatives of the (finitely many) generators of $\Gamma$. If $\tau$ 
belongs to $]0,1[$ and $\Gamma$ consists of $C^{1+\tau}$ diffeomorphisms, 
the same is true provided that the sum of the $\tau$-powers of the lengths 
of the corresponding 
images of $I$ is finite (this last condition does not follow from 
the disjointness of these intervals~!): see for instance \cite{DKN}, 
Lemma 2.2. It is not difficult to prove a similar statement for the 
case of different regularities, and this is precisely the content of the 
following lemma. However, to the difference of \cite{DKN}, here 
we will deal with {\em finite} sequences of compositions by a 
technical reason which will be clear at the end of the 
next section.

\vspace{0.1cm}

\begin{lem} {\em Let $\Gamma$ be a semigroup of (orientation 
preserving) diffeomorphisms 
of the circle or the closed interval which is generated by finitely 
many elements $g_k$, $k \in \{1,\ldots,l\}$, which are respectively of 
class $C^{1+\tau_k}$, where $\tau_k \! \in ]0,1]$. Let $C_k$ denote 
the $\tau_k$-H\"older constant of the function $\log(g_k')$, and let 
\esp $C=\max \{C_1,\ldots,C_l\}$ \esp and \esp 
$\tau = \max \{\tau_1,\ldots,\tau_l\}$. \esp 
Given $n_0 \in \mathbb{N}$, for each $n \leq n_0$ let us chose 
$k_n \in \{1,\ldots,l\}$, and for a fixed  interval 
$I$ let $S>0$ be a constant such that 
\begin{equation}
S \geq \sum_{n = 0}^{n_0-1} \big| g_{k_n} \cdots g_{k_1} (I) \big|^{\tau_{k_{n+1}}}.
\label{sum-eme}
\end{equation}
If $n \leq n_0$ is such that $g_{k_n} \cdots g_{k_1} (I)$ does not intersect $I$ but 
is contained in the $L$-neighborhood of $I$, where $L := |I| / 2 \exp(2^{\tau}CS)$, 
then $g_{k_n} \cdots g_{k_1}$ has a hyperbolic fixed point.} 
\label{lema-cle}
\end{lem}

\noindent{\bf Proof.} Let $J = [a,b]$ be the (closed) $2L$-neighborhood of $I$, and 
let $I'$ (resp. $I''$) the connected component of $J \setminus I$ to the right (resp. 
to the left) of $I$. We will prove by induction on $j \!\in\! \{0,\ldots,n_0\}$ 
that the following two conditions are satisfied:

\vspace{0.1cm}

\noindent{$(\mathrm{i})_j \hspace{0.25cm} |g_{k_{j}} \cdots g_{k_{1}}(I')| \leq  
|g_{k_{j}} \cdots g_{k_{1}}(I)|$,}

\vspace{0.1cm}

\noindent{$(\mathrm{ii})_j \hspace{0.25cm} \sup_{ \{ x,y \} \subset I \cup I'} 
\frac{(g_{k_{j}} \cdots g_{k_{1}})'(x)}{(g_{k_{j}} \cdots g_{k_{1}})'(y)} 
\leq \exp(2^{\tau} \thinspace C S )$.}

\vspace{0.1cm}

Condition $(\mathrm{ii})_0$ is trivially satisfied, whereas condition $(\mathrm{i})_0$ 
is satisfied since $|I'| \!=\! 2L \!\leq\! |I|$. Assume that $(\mathrm{i})_i$ and 
$(\mathrm{ii})_i$ hold for each $i \in \{0,\ldots,j-1\}$. Then 
for every $x,y$ in $I \cup I'$ we have
\begin{eqnarray*}
\left|\log\left(\frac{(g_{k_{j}} \cdots g_{k_{1}})'(x)}
{(g_{k_{j}} \cdots g_{k_{1}})'(y)}\right)\right| 
&\leq& \sum_{i=0}^{j-1} \big| \log(g_{k_{i+1}}'(g_{k_{i}} \cdots g_{k_{1}}(x))) - 
\log(g_{k_{i+1}}'(g_{k_{i}} \cdots g_{k_{1}}(y))) \big|\\
&\leq& \esp \sum_{i=0}^{j-1} C_{k_{i+1}} \big| g_{k_{i}} \cdots g_{k_{1}}(x) - 
g_{k_{i}} \cdots g_{k_{1}}(y) \big|^{\tau_{k_{i+1}}}\\
&\leq& C \esp \esp \sum_{i=0}^{j-1} \big( |g_{k_{i}} \cdots g_{k_{1}}(I)| 
+ |g_{k_{i}} \cdots g_{k_{1}}(I')| \big)^{\tau_{k_{i+1}}}\\
&\leq& C \esp \esp 2^{\tau} \sum_{i=0}^{j-1} |g_{k_i} \cdots g_{k_1}(I)|^{\tau_{k_{i+1}}}\\
&\leq& C \esp \esp 2^{\tau} S.
\end{eqnarray*}
This shows $(\mathrm{ii})_j$. To verify $(\mathrm{i})_j$ first 
note that there must exist $x \in I$ and $y \in I'$ such that 
$$|g_{k_{j}} \cdots g_{k_{1}} (I)| = |I| \cdot (g_{k_{j}} \cdots g_{k_{1}})'(x) 
\qquad \mbox { and } \qquad |g_{k_{j}} \cdots g_{k_{1}}(I')| 
= |I'| \cdot (g_{k_{j}} \cdots g_{k_{1}})'(y).$$
Therefore, by $(\mathrm{ii})_j$,
$$\frac{|g_{k_{j}} \cdots g_{k_{1}}(I')|}{|g_{k_{j}} \cdots g_{k_{1}}(I)|}= 
\frac{(g_{k_{j}} \cdots g_{k_{1}})'(x)}{(g_{k_{j}} \cdots g_{k_{1}})'(y)} 
\cdot \frac{|I'|}{|I|} \leq \exp(2^{\tau} C S) \frac{|I'|}{|I|} \leq 1,$$
which proves $(\mathrm{i})_j$. Obviously, similar arguments show that 
$(\mathrm{i})_j$ and $(\mathrm{ii})_j$ also hold for every 
$j\!\!~\in~\!\!\{0,\ldots,n_0\}$ when we replace $I'$ by $I''$.

Now for simplicity let us denote $h_j = g_{k_j} \cdots g_{k_1}$. Assume 
that $h_n(I)$ is contained in the $L$-neighborhood of the interval 
$I$ (see Figure 1). Then property $(\mathrm{i})_n$ gives $h_n(J) \subset J$, 
and this already implies that $h_n$ has a fixed point $x$ in $J$. 
(The reader will see that the existence of this fixed point 
together with the fact that $h_n \neq id$ is the only information 
that we will retain for the proof of Theorem A.) 

To conclude we would like to show that the fixed point $x$ is hyperbolic. To do this just 
note that, if $h_n(I)$ does not intersect $I$, then there exists $y \in I$ such that 
$$h_n'(y) = \frac{|h_n(I)|}{|I|} \leq \frac{L}{|I|}.$$
Therefore, by $(\mathrm{ii})_n$,
$$h_n'(x) \leq h_n'(y) \exp(2^{\tau} C S) \leq 
\frac{L \exp(2^{\tau}CS)}{|I|} \leq \frac{1}{2},$$
and this finishes the proof. $\hfill\square$ 

\vspace{0.35cm}

\beginpicture

\setcoordinatesystem units <1cm,1cm>

\putrule from -5.5 0 to 5.5 0 

\putrule from -2.8 -0.015 to 2.8 -0.015

\putrule from -2.8 0.015 to 2.8 0.015

\putrule from 3 -0.015 to 3.5 -0.015 

\putrule from 3 0.015 to 3.5 0.015

\circulararc -45 degrees from 0 0.5 center at 1.6 -3.3

\circulararc 132 degrees from 3.8 -0.5 center at 4.25 -0.3 

\circulararc -55 degrees from 2.2 -0.5 center at -1.25 6

\plot 3.8 -0.5 
3.9 -0.56 /

\plot 3.8 -0.5 
3.8 -0.6 /

\plot 2.2 -0.5 
2.12 -0.63 /

\plot 2.2 -0.5 
2.03 -0.53 /

\plot 3.14 0.52 
3 0.52 /

\plot 3.14 0.52 
3.05 0.62 /

\put{$\Big|$} at -4.8 0  
\put{$\Big|$} at 4.8 0 
\put{$|$} at 3.7 0  
\put{$|$} at 2.5 0
\put{$\Big($} at -2.8 0
\put{$\Big)$} at 2.8 0 
\put{$($} at 3 0
\put{$)$} at 3.5 0 
\put{$h_n$} at -3.8 -1.2
\put{$h_n$} at 4.4 -1.2
\put{$h_n$} at 1.3 1.1
\put{$a$} at -5.051 0.23 
\put{$b$} at 5.051 0.23

\plot 3.615 1 
3.615 0.3 /

\plot 3.615 0.3
3.55 0.45 /

\plot 3.615 0.3 
3.68 0.45 /

\put{Figure 1} at 0 -1.8 
\put{} at -8.1 0

\small
\put{$h_n(I)$} at 3.25 -0.6
\put{$I$} at 0 -0.6 
\put{$\bullet$} at 3.62 0 
\put{hyperbolic} at 3.615 1.6
\put{fixed point} at 3.615 1.25

\endpicture


\vspace{0.45cm}

\section{Proof of Theorem A}

\hspace{0.5cm} Recall the following well known argument (see for instance \cite{ghys}, 
Proposition 6.17, or \cite{imca}, Lemma 4.1.4). If $f_1,\ldots,f_d$ are commuting circle 
homeomorphisms, then there is a common invariant probability measure $\mu$ on $\clo$. 
Moreover, if the rotation number of at least one of them is irrational, 
then there is no finite orbit for the group action, 
and the measure $\mu$ has no atom. Therefore, the distribution function 
$$F_{\mu}:S^1 \to \mathbb{R}/\mathbb{Z}, \qquad  F_{\mu}(x):=\mu([0,x[),$$
gives a (simultaneous) semiconjugacy between the maps $f_1,\dots,f_d$ 
and the rotations corresponding to their rotation numbers. Thus, for the proof of  
Theorem A we have to show that this semiconjugacy is in fact a conjugacy, and our 
strategy for proving this (under the hypothesis of the Theorem) is the classical 
one and goes back to Schwartz \cite{Sc}. Indeed, in the contrary case the support 
of $\mu$ would be a (minimal) invariant Cantor set, and the connected components 
of its complement would correspond to the maximal wandering open intervals. 
Fixing one of these intervals, say $I$, we will search for 
a sequence of compositions $h_n = f_{k_n} \cdots f_{k_1}$ 
satisfying the hypothesis of Lemma \ref{lema-cle}. This will allow us 
to conclude that some $h_n$ has a (hyperbolic) fixed point, thus implying 
that its rotation number is equal to zero. 
However, this is in contradiction to the fact that the rotation 
numbers of the $f_k$'s are independent over the rationals 
(it is easy to verify that the rotation number restricted to any 
group of circle homeomorphisms which preserves a probability 
measure on $\clo$ is a group homomorphism: see again 
\cite{ghys} or \cite{imca}).

In order to ensure the existence of the sequence $(h_n)$ the main idea of 
\cite{DKN} was to endow the space of all (infinite) sequences of compositions with 
a natural probability measure, and then to prove that the ``generic ones" satisfy 
many nice properties as for instance the convergence of the sum (\ref{sum-eme}) 
as $n_0$ goes to infinity. It seems that such a probabilistic argument cannot be 
applied to the case of different regularities, and we will need to introduce a new 
argument which is somehow more deterministic, since it gives partial information 
on the sequence that we find. For simplicity we will first deal with the case $d\!=\!2$.


\subsection{The case $d=2$}

\hspace{0.5cm} Although not explicitly stated in \cite{DKN}, the main 
probabilistic argument for the proof of the Generalized Denjoy Theorem 
therein is not a dynamical issue, but it is just a statement concerning 
the finiteness of the sum of the $\tau$-powers of some positive real numbers. 
To be more concrete (at least in the case $d=2$ and when $\tau > 1/2$), if 
$(\ell_{i,j})$ is a double-indexed sequence of positive numbers with finite 
total sum (where $i$ and $j$ are non negative integers), then with 
respect to some natural probability distribution on the space of 
infinite paths $(i(n),j(n))_{n \geq 0}$ satisfying 
$i(0)=j(0)=0$, $i(n+1) \geq i(n)$, $j(n+1) \geq j(n)$ and $i(n+1)+j(n+1) 
= 1 + i(n) + j(n)$, one has almost everywhere the convergence of the sum 
$$\sum_{n \geq 0} \ell_{i(n),j(n)}^{\tau}.$$ 
The first goal of this section is to prove the existence of paths sharing 
a similar property in the case of different exponents $\tau_1,\tau_2$ in 
$]0,1[$ (with $\tau_1 + \tau_2 > 1$). A substantial difference here is 
that we will construct our sequence by concatenating infinitely many 
finite paths, and each one of these paths will be chosen among finitely 
many ones. To do this we begin with the following elementary lemma.

\vspace{0.1cm}

\begin{lem} {\em Let $\ell_{i,j}$ be positive real numbers, where $i \in \{1,\ldots,m\}$ 
and $j \in \{1,\ldots,n\}$. Assume that the total sum of the $\ell_{i.j}$'s is less 
than or equal to $1$. If $\tau$ belongs to $]0,1[$, then there exists 
$k \in \{1,\ldots,n\}$ such that}
$$\sum_{i = 1}^{m} \ell_{i,k}^{\tau} \leq \frac{m^{1-\tau}}{n^{\tau}}.$$
\label{basico}
\end{lem}

\noindent{\bf Proof.} We will show that the mean value of the function
\esp \esp $k \mapsto \sum_{i = 1}^{m} \ell_{i,k}^{\tau}$ \esp \esp 
is less than or equal to $m^{1-\tau}/n^{\tau}$, from where the 
claim of the lemma follows immediately. To do this first note that, by 
H\"older's inequality, for each fixed $k \in \{1,\ldots,n\}$ one has 
$$\sum_{i = 1}^{m} \ell_{i,k}^{\tau} = 
\left\langle (\ell_{i,k}^{\tau})_{i=1}^{m}, 
(1)_{i=1}^{m} \right\rangle \leq 
\left\| (\ell_{i,k}^{\tau})_{i=1}^{m} \right\|_{1/\tau} \cdot 
\left\| (1)_{i=1}^{m} \right\|_{1/(1-\tau)} = 
\left(\sum_{i=1}^{m} \ell_{i,k} \right)^{\tau} m^{1-\tau}.$$
Thus, by using H\"older's inequality again one obtains
\begin{eqnarray*}
\frac{1}{n} \sum_{k=1}^{n} \left( \sum_{i = 1}^{m} \ell_{i,k}^{\tau} \right) 
&=& \frac{m^{1-\tau}}{n} 
\left\langle \left( \Big(\sum_{k=1}^{n} \ell_{i,k}\Big)^{\tau} \right)_{k=1}^{n} , 
\left( 1 \right)_{k=1}^{n} \right\rangle\\
&\leq& \frac{m^{1-\tau}}{n} 
\left\| \left( \Big(\sum_{k=1}^{n} \ell_{i,k} \Big)^{\tau}\right)_{k=1}^{n} 
\right\|_{1/\tau} \cdot \left\| \left( 1 \right)_{k=1}^{n} \right\|_{1/(1-\tau)}\\
&=& \frac{m^{1-\tau}}{n} \left( \sum_{k=1}^{n} \sum_{i=1}^{m} \ell_{i,k} \right)^{\tau} 
n^{1-\tau}\\ 
&\leq& \frac{m^{1-\tau}}{n^{\tau}},
\end{eqnarray*}
which finishes the proof. $\hfill\square$

\vspace{0.38cm}

Now we explain the main idea of our construction. Let us assume that the total sum of 
the double-indexed sequence of positive numbers $\ell_{i,j}$ is $\leq 1$, and suppose 
that the numbers $\tau_1 \! \in ]0,1[$ and $\tau_2 \! \in ]0,1[$ such that 
$\tau_1 +\tau_2 > 1$ are fixed. Denoting by $[[a,b]]$ the set of integers 
between $a$ and $b$ (with $a$ and $b$ included when they 
are in $\mathbb{Z}$), let us consider any sequence 
of rectangles \esp $R_m \subset \mathbb{N}_0 \times \mathbb{N}_0$ \esp 
such that 
\esp $R_0 = \{(0,0)\}$, \esp $R_{2m+1} = [[i_{m},i_{m+1}]] 
\times [[j_{m},j_{m+2}]]$ \esp and \esp 
$R_{2m+2}~=~[[i_{m},i_{m+2}]] \times [[j_{m+1},j_{m+2}]]$, 
where $(i_m)_{m \geq 1}$ \esp and $(j_m)_{m \geq 1}$ are 
strictly increasing sequences of non negative integers 
numbers satisfying $i_0 \!=\! i_1 \!=\! 0$ and $j_0 \!=\! j_1 \!=\! 0$  
(see Figure 2). Denoting by $X_m$ and $Y_m$ respectively the number of 
points on the horizontal and vertical sides of each $R_m$, a direct application 
of Lemma \ref{basico} gives us, for $\varepsilon := 1 - \tau_1 - \tau_2 > 0$ 
and each $m \geq 0$:

\vspace{0.1cm}

\noindent -- an integer $r(2m+1) \in [[i_{m},i_{m+1}]]$ such that 
$$\sum_{j=j_{m}}^{j_{m+2}} \ell_{r(2m+1),j}^{\tau_2} \leq 
\frac{Y_{2m+1}^{1 - \tau_2}}{X_{2m+1}^{\tau_2}} = 
\frac{Y_{2m+1}^{\tau_1}}{X_{2m+1}^{\tau_2}} \cdot Y_{2m+1}^{-\varepsilon},$$ 

\vspace{0.1cm}

\noindent -- an integer $r(2m+2) \in [[j_{m+1},j_{m+2}]]$ such that 
$$\sum_{i=i_{m}}^{i_{m+2}} \ell_{i,r(2m+2)}^{\tau_1} \leq 
\frac{X_{2m+2}^{1-\tau_1}}{Y_{2m+2}^{\tau_1}} = 
\frac{X_{2m+2}^{\tau_2}}{Y_{2m+2}^{\tau_1}} \cdot X_{2m+2}^{-\varepsilon}.$$ 

Starting from the origin and following the corresponding horizontal and 
vertical lines, we find an infinite path $(i(n),j(n))_{n \geq 0}$ satisfying 
$$i(0) = j(0) = 0, \quad i(n+1) \geq i(n), \quad j(n+1) \geq j(n), 
\quad i(n+1)+j(n+1) = 1 + i(n) + j(n),$$ 
and such that the sum 
\begin{equation}
\sum_{n \geq 0} \ell_{i(n),j(n)}^{\tau_{\alpha(n)}}
\label{acontrolar}
\end{equation} 
is bounded by 
\begin{equation}
\sum_{m \geq 0} \left[ \frac{Y_{2m+1}^{\tau_1}}{X_{2m+1}^{\tau_2}} 
\cdot Y_{2m+1}^{-\varepsilon} + \frac{X_{2m+2}^{\tau_2}}{Y_{2m+2}^{\tau_1}} 
\cdot X_{2m+2}^{-\varepsilon} \right],
\label{control}
\end{equation}
where $\alpha(n) := 1$ if \esp $|i(n+1) - i(n)| = 1$ \esp and 
$\alpha(n) := 2$ if \esp $|j(n+1) - j(n)| = 1$.

\vspace{0.86cm}


\beginpicture

\setcoordinatesystem units <1.5cm,1.5cm>

\putrule from 0 0 to 0.5 0 

\putrule from 0 0.49 to 0.5 0.49 
\putrule from 0 0.51 to 0.5 0.51

\putrule from 0.01 0 to 0.01 0.5  
\putrule from -0.01 0 to -0.01 0.5

\putrule from 0 0.5 to 1.5 0.5  
\putrule from 0 1.2 to 4 1.2 
\putrule from 0 0 to 0 1.2    
\putrule from 0.5 0 to 0.5 2.5   
\putrule from 1.5 0.5 to 1.5 5.6  
\putrule from 4 1.2 to 4 6.01  
\putrule from 0.5 2.5 to 8 2.5 
\putrule from 1.5 5.6 to 8 5.6 
\putrule from 8 2.5 to 8 6.01   

\putrule from 0.5 0 to 0.5 1.2
\putrule from 0.49 0 to 0.49 1.2 
\putrule from 0.51 0 to 0.51 1.2 

\putrule from 0 0.8 to 1.5 0.8
\putrule from 0 0.81 to 1.5 0.81
\putrule from 0 0.82 to 1.5 0.82 

\putrule from 0.8 0.5 to 0.8 2.5 
\putrule from 0.81 0.5 to 0.81 2.5 
\putrule from 0.82 0.5 to 0.82 2.5 

\putrule from 0.5 2 to 4 2 
\putrule from 0.5 2.01 to 4 2.01
\putrule from 0.5 2.02 to 4 2.02 

\putrule from 3 1.2 to 3 5.6 
\putrule from 3.01 1.2 to 3.01 5.6 
\putrule from 3.02 1.2 to 3.02 5.6 

\putrule from 1.5 4.5 to 8 4.5 
\putrule from 1.5 4.51 to 8 4.51  
\putrule from 1.5 4.52 to 8 4.52  

\putrule from 5.7965 2.5 to 5.7965 6 
\putrule from 5.8065 2.5 to 5.8065 6
\putrule from 5.8165 2.5 to 5.8165 6

\put{Figure 2} at 4 -0.75

\setdots

\putrule from 0.5 2.5 to 0.5 5.6 
\putrule from 0.5 0 to 8 0 
\putrule from 1.5 0 to 1.5 0.5 
\putrule from 4 0 to 4 1.2 
\putrule from 8 0 to 8 2.5 

\putrule from 0 2.5 to 0.5 2.5 
\putrule from 0 5.6 to 1.5 5.6 

\putrule from 0 1.2 to 0 5.6 
\putrule from 0.81 2.5 to 0.81 5.6 

\putrule from 1.5 0.81 to 8 0.81  
\putrule from 4 2.01 to 8 2.01  

\putrule from 1.5 0.5 to 8 0.5 

\begin{tiny}
\put{$i_0\!=\!i_1$} at 0 -0.23 
\put{$i_2$} at 0.5 -0.23 
\put{$i_3$} at 1.5 -0.23 
\put{$i_4$} at 4 -0.23  
\put{$i_5$} at 8 -0.23  
\put{$j_0\!=\!j_1$} at -0.35 0 
\put{$j_2$} at -0.35 0.5  
\put{$j_3$} at -0.35 1.2 
\put{$j_4$} at -0.35 2.5 
\put{$j_5$} at -0.35 5.6 

\put{$r(1)$} at -0.1 5.9 
\put{$r(0)$} at 8.3 0 
\put{$r(2)$} at 8.3 0.5 
\put{$r(3)$} at 0.45 5.9 
\put{$r(4)$} at 8.3 0.81 
\put{$r(5)$} at 0.85 5.9  
\put{$r(6)$} at 8.3 2.01   
\put{$r(7)$} at 3.01 5.9 
\put{$r(8)$} at 8.3 4.51 

\put{$\bullet$} at -0.05 0  
\put{$\bullet$} at -0.05 0.1  
\put{$\bullet$} at -0.05 0.2  
\put{$\bullet$} at -0.05 0.3  
\put{$\bullet$} at -0.05 0.4  
\put{$\bullet$} at -0.05 0.5 
  
\put{$\bullet$} at 0.445 0.5 
\put{$\bullet$} at 0.445 0.6  
\put{$\bullet$} at 0.445 0.7 
\put{$\bullet$} at 0.445 0.8 

\put{$\bullet$} at -0.05 0.5
\put{$\bullet$} at 0.05 0.5  
\put{$\bullet$} at 0.15 0.5  
\put{$\bullet$} at 0.25 0.5  
\put{$\bullet$} at 0.35 0.5  

\put{$\bullet$} at 0.45 0.81  
\put{$\bullet$} at 0.55 0.81  
\put{$\bullet$} at 0.65 0.81  
\put{$\bullet$} at 0.75 0.81   
   
\put{$\bullet$} at 0.75 0.91  
\put{$\bullet$} at 0.75 1.01  
\put{$\bullet$} at 0.75 1.11  
\put{$\bullet$} at 0.75 1.21  
\put{$\bullet$} at 0.75 1.31  
\put{$\bullet$} at 0.75 1.41  
\put{$\bullet$} at 0.75 1.51  
\put{$\bullet$} at 0.75 1.61  
\put{$\bullet$} at 0.75 1.71  
\put{$\bullet$} at 0.75 1.81  
\put{$\bullet$} at 0.75 1.91 
\put{$\bullet$} at 0.75 2.01    

\put{$\bullet$} at 0.85 2.01 
\put{$\bullet$} at 0.95 2.01 
\put{$\bullet$} at 1.05 2.01  
\put{$\bullet$} at 1.15 2.01  
\put{$\bullet$} at 1.25 2.01   
\put{$\bullet$} at 1.35 2.01   
\put{$\bullet$} at 1.45 2.01  
\put{$\bullet$} at 1.55 2.01  
\put{$\bullet$} at 1.65 2.01  
\put{$\bullet$} at 1.75 2.01  
\put{$\bullet$} at 1.85 2.01 
\put{$\bullet$} at 1.95 2.01  
\put{$\bullet$} at 2.05 2.01  
\put{$\bullet$} at 2.15 2.01  
\put{$\bullet$} at 2.25 2.01   
\put{$\bullet$} at 2.35 2.01   
\put{$\bullet$} at 2.45 2.01  
\put{$\bullet$} at 2.55 2.01  
\put{$\bullet$} at 2.65 2.01  
\put{$\bullet$} at 2.75 2.01  
\put{$\bullet$} at 2.85 2.01 
\put{$\bullet$} at 2.95 2.01    

\put{$\bullet$} at 2.95 2.11 
\put{$\bullet$} at 2.95 2.21 
\put{$\bullet$} at 2.95 2.31  
\put{$\bullet$} at 2.95 2.41  
\put{$\bullet$} at 2.95 2.51   
\put{$\bullet$} at 2.95 2.61   
\put{$\bullet$} at 2.95 2.71  
\put{$\bullet$} at 2.95 2.81  
\put{$\bullet$} at 2.95 2.91  
\put{$\bullet$} at 2.95 3.01  
\put{$\bullet$} at 2.95 3.11 
\put{$\bullet$} at 2.95 3.21  
\put{$\bullet$} at 2.95 3.31  
\put{$\bullet$} at 2.95 3.41  
\put{$\bullet$} at 2.95 3.51   
\put{$\bullet$} at 2.95 3.61   
\put{$\bullet$} at 2.95 3.71  
\put{$\bullet$} at 2.95 3.81  
\put{$\bullet$} at 2.95 3.91  
\put{$\bullet$} at 2.95 4.01  
\put{$\bullet$} at 2.95 4.11 
\put{$\bullet$} at 2.95 4.21    
\put{$\bullet$} at 2.95 4.31 
\put{$\bullet$} at 2.95 4.41 
\put{$\bullet$} at 2.95 4.51 

\put{$\bullet$} at 3.05 4.51  
\put{$\bullet$} at 3.15 4.51  
\put{$\bullet$} at 3.25 4.51  
\put{$\bullet$} at 3.35 4.51 
\put{$\bullet$} at 3.45 4.51    
\put{$\bullet$} at 3.55 4.51 
\put{$\bullet$} at 3.65 4.51 
\put{$\bullet$} at 3.75 4.51 
\put{$\bullet$} at 3.85 4.51  
\put{$\bullet$} at 3.95 4.51  
\put{$\bullet$} at 4.05 4.51  
\put{$\bullet$} at 4.15 4.51 
\put{$\bullet$} at 4.25 4.51    
\put{$\bullet$} at 4.35 4.51 
\put{$\bullet$} at 4.45 4.51 
\put{$\bullet$} at 4.55 4.51 
\put{$\bullet$} at 4.65 4.51    
\put{$\bullet$} at 4.75 4.51 
\put{$\bullet$} at 4.85 4.51 
\put{$\bullet$} at 4.95 4.51 
\put{$\bullet$} at 5.05 4.51  
\put{$\bullet$} at 5.15 4.51  
\put{$\bullet$} at 5.25 4.51  
\put{$\bullet$} at 5.35 4.51 
\put{$\bullet$} at 5.45 4.51    
\put{$\bullet$} at 5.55 4.51 
\put{$\bullet$} at 5.65 4.51 
\put{$\bullet$} at 5.75 4.51 

\put{$\bullet$} at 5.75 4.61 
\put{$\bullet$} at 5.75 4.71 
\put{$\bullet$} at 5.75 4.81 
\put{$\bullet$} at 5.75 4.91 
\put{$\bullet$} at 5.75 5.01  
\put{$\bullet$} at 5.75 5.11  
\put{$\bullet$} at 5.75 5.21  
\put{$\bullet$} at 5.75 5.31  
\put{$\bullet$} at 5.75 5.41  
\put{$\bullet$} at 5.75 5.51  
\put{$\bullet$} at 5.75 5.61  
\put{$\bullet$} at 5.75 5.71  
\put{$\bullet$} at 5.75 5.81  
\put{$\bullet$} at 5.75 5.91  
\put{$\bullet$} at 5.75 6.01  

\put{$R_2$} at 0.3 0.38 
\put{$R_3$} at 0.3 1.1 
\put{$R_5$} at 1.3 2.4 
\put{$R_6$} at 3.8 2.4 
\put{$R_7$} at 3.8 5.5 
\put{$R_8$} at 7.8 5.5 
\put{$R_4$} at 1.3 1.1 
\end{tiny} 

\put{} at -1.4 2.5 

\endpicture


\vspace{0.7cm}

Now let us consider any choice such that $i_m = [4^{m\tau_1}]$ and 
$j_m = [4^{m\tau_2}]$ for $m$ large enough. Writing $a_m \simeq b_m$ 
when $(a_m)$ and $(b_m)$ are sequences of positive numbers such that 
$(a_m / b_m)$ remains bounded and away from zero, for such a choice we 
have $X_m \simeq 2^{m \tau_1}$ and $Y_m \simeq 2^{m \tau_2}$. Thus, 
$$\frac{X_{m}^{\tau_2}}{Y_{m}^{\tau_1}} \simeq
\frac{(2^{m \tau_1 })^{\tau_2}}{(2^{m \tau_2 })^{\tau_1}} = 1,$$
and therefore there exists $C>0$ such that, for each $m \geq 0$,
$$\frac{1}{C} \leq \frac{X_{m}^{\tau_2}}{Y_{m}^{\tau_1}} \leq C.$$ 
This implies that the sum in (\ref{control}) is bounded by 
\begin{eqnarray}
S := C \left( \sum_{m \geq 0} 
\Big[ \Big( \frac{1}{4^{m\tau_2}} \Big)^{\varepsilon} 
+ \Big( \frac{1}{4^{m\tau_1}} \Big)^{\varepsilon} \Big] \right) 
= C \left( \frac{4^{\tau_2 \varepsilon}}{4^{\tau_2 \varepsilon} - 1} + 
\frac{4^{\tau_1 \varepsilon}}{4^{\tau_1 \varepsilon} - 1} \right),
\label{def-ese}
\end{eqnarray}
and so the value of the sum (\ref{acontrolar}) 
is finite (and also bounded by $S$).

\vspace{0.25cm}

We can now proceed to the proof of Theorem A in the case $d \!=\! 2$. Assume by 
contradiction that $f_k, k \!\in\! \{1,2\}$, are respectively $C^{1+\tau_k}$ 
commuting circle diffeomorphisms which are not simultaneously conjugate to 
rotations and which have rotation numbers independent over the rationals. Let 
$I$ be a connected component of the complement of the invariant minimal Cantor 
set for the group action, and let $\ell_{i,j} = |f_1^if_2^j(I)|$. We obviously 
have $\sum_{i,j} \ell_{i,j} \leq 1$, and so we can apply all our previous 
discussion to this sequence. In particular, there exists an infinite path 
$(i(n),j(n))$ starting at the origin and such that the sum  
$$\sum_{n \geq 0} \ell_{i(n),j(n)}^{\tau_{\alpha(n)}}$$
is bounded by the number $S > 0$ defined by (\ref{def-ese}). 
If for $n \geq 1$ we let $k_n = \alpha(n-1) \in \{ 1,2\}$, then 
we obtain a sequence of compositions $h_n = f_{k_n} \cdots f_{k_1}$ 
such that the preceding sum coincides term by term with
$$\sum_{n \geq 0} |f_{k_n} \cdots f_{k_1} (I)|^{\tau_{k_{n+1}}}.$$
Thus, in order to apply Lemma \ref{lema-cle} to get a contradiction, we 
just need to verify that, for some $n \geq 1$, the hypothesis 
that $h_n(I) = f_{k_n} \cdots f_{k_1} (I)$ is contained in 
the $L$-neighborhood of $I$ is satisfied (where $L:=|I
|/2\exp(2^{\tau} C S)$, 
$\tau := \max \{\tau_1,\tau_2\}$, and $C := \max \{C_1,\ldots,C_d\}$, with 
$C_k$ being the $\tau_k$-H\"older constant for the function $\log(f_k')$). 

To to this first note that, if we collapse all the connected components of 
the complement of the minimal invariant Cantor set, then we obtain a 
topological circle $\hat{\mathrm{S}}^1$ on which the original diffeomorphisms 
induce naturally minimal homeomorphisms $\hat{f}_1$ and $\hat{f}_2$ 
which are simultaneously conjugate to rotations. Moreover, 
the $L$-neighborhood of $I$ becomes a non degenerate interval 
$\hat{U}$; thus, there exists $N \in \mathbb{N}$ such that the intervals 
$\hat{f}_1^{-1}(\hat{U}),\ldots,\hat{f}_1^{-N}(\hat{U})$, as well as 
$\hat{f}_2^{-1}(\hat{U}),\ldots,\hat{f}_2^{-N}(\hat{U})$, cover the 
circle $\hat{\mathrm{S}}^1$. This easily implies that for any image 
$I_0$ of $I$ by some element of the semigroup generated by $f_1$ 
and $f_2$ there exists $k$ and $k'$ in $\{1,\ldots,N\}$ such that 
$f_1^k(I_0)$ and $f_2^{k'}(I_0)$ are contained in the $L$-neighborhood 
of $I$. Now it is easy to see that, for the sequence of compositions that 
we found, for every $\bar{N} \in \mathbb{N}$ there exists some integer 
$r \in \mathbb{N}$ such that $k_r = k_{r+1} = \ldots = k_{r+\bar{N}}$. 
For $\bar{N} = N$ this obviously implies that at least one of the intervals 
$h_{r+1}(I),\ldots,h_{r+N}(I)$ is contained in the $L$-neighborhood of $I$, 
thus finishing the proof.

\vspace{1cm}


\beginpicture

\setcoordinatesystem units <1.5cm,1.5cm>

\putrule from 0 0 to 8 0   
\putrule from 0 0.5 to 1.5 0.5  
\putrule from 0 2.2 to 3.5 2.2       
\putrule from 0 4.5 to 8 4.5 

\putrule from 0 0 to 0 4.5    
\putrule from 0.5 0 to 0.5 0.5   
\putrule from 1.5 0 to 1.5 2.2   
\putrule from 3.5 0 to 3.5 4.5 
\putrule from 8 0 to 8 4.5 

\put{Figure 3} at 4 -0.58

\putrule from -0.01 0 to -0.01 0.5 
\putrule from 0.01 0 to 0.01 0.5 

\putrule from 0 0.2 to 1.5 0.2
\putrule from 0 0.21 to 1.5 0.21
\putrule from 0 0.22 to 1.5 0.22 

\putrule from 0.8 0 to 0.8 2.2 
\putrule from 0.81 0 to 0.81 2.2 
\putrule from 0.82 0 to 0.82 2.2 

\putrule from 0 1 to 3.5 1 
\putrule from 0 1.01 to 3.5 1.01
\putrule from 0 1.02 to 3.5 1.02 

\putrule from 2 0 to 2 4.5 
\putrule from 2.01 0 to 2.01 4.5 
\putrule from 2.02 0 to 2.02 4.5 

\putrule from 0 0.7 to 8 0.7 
\putrule from 0 0.71 to 8 0.71  
\putrule from 0 0.72 to 8 0.72  

\begin{tiny}
\put{$\bullet$} at -0.055 0   
\put{$\bullet$} at -0.055 0.1  
\put{$\bullet$} at -0.055 0.2    
\put{$\bullet$} at 0.045 0.2  
\put{$\bullet$} at 0.145 0.2   
\put{$\bullet$} at 0.245 0.2 
\put{$\bullet$} at 0.345 0.2  
\put{$\bullet$} at 0.445 0.2   
\put{$\bullet$} at 0.545 0.2  
\put{$\bullet$} at 0.645 0.2   
\put{$\bullet$} at 0.745 0.2    
 
\put{$\bullet$} at 0.75 0.3   
\put{$\bullet$} at 0.75 0.4    
\put{$\bullet$} at 0.75 0.5   
\put{$\bullet$} at 0.75 0.6
\put{$\bullet$} at 0.75 0.7  
\put{$\bullet$} at 0.75 0.8   
\put{$\bullet$} at 0.75 0.9    
\put{$\bullet$} at 0.75 1   

\put{$\bullet$} at 0.85 1
\put{$\bullet$} at 0.95 1    
\put{$\bullet$} at 1.05 1 
\put{$\bullet$} at 1.15 1 
\put{$\bullet$} at 1.25 1  
\put{$\bullet$} at 1.35 1  
\put{$\bullet$} at 1.45 1   
\put{$\bullet$} at 1.55 1   
\put{$\bullet$} at 1.65 1  
\put{$\bullet$} at 1.75 1  
\put{$\bullet$} at 1.85 1  
\put{$\bullet$} at 1.95 1
  
\put{$\bullet$} at 1.95 0.9 
\put{$\bullet$} at 1.95 0.8  
\put{$\bullet$} at 1.95 0.7

\put{$\bullet$} at 1.95 0.7
\put{$\bullet$} at 2.05 0.7
\put{$\bullet$} at 2.15 0.7
\put{$\bullet$} at 2.25 0.7  
\put{$\bullet$} at 2.35 0.7 
\put{$\bullet$} at 2.45 0.7
\put{$\bullet$} at 2.55 0.7
\put{$\bullet$} at 2.65 0.7
\put{$\bullet$} at 2.75 0.7  
\put{$\bullet$} at 2.85 0.7
\put{$\bullet$} at 2.95 0.7
\put{$\bullet$} at 3.05 0.7
\put{$\bullet$} at 3.15 0.7
\put{$\bullet$} at 3.25 0.7  
\put{$\bullet$} at 3.35 0.7 
\put{$\bullet$} at 3.45 0.7
\put{$\bullet$} at 3.55 0.7
\put{$\bullet$} at 3.65 0.7
\put{$\bullet$} at 3.75 0.7  
\put{$\bullet$} at 3.85 0.7
\put{$\bullet$} at 3.95 0.7
\put{$\bullet$} at 4.05 0.7
\put{$\bullet$} at 4.15 0.7
\put{$\bullet$} at 4.25 0.7  
\put{$\bullet$} at 4.35 0.7 
\put{$\bullet$} at 4.45 0.7
\put{$\bullet$} at 4.55 0.7
\put{$\bullet$} at 4.65 0.7
\put{$\bullet$} at 4.75 0.7  
\put{$\bullet$} at 4.85 0.7
\put{$\bullet$} at 4.95 0.7
\put{$\bullet$} at 5.05 0.7
\put{$\bullet$} at 5.15 0.7
\put{$\bullet$} at 5.25 0.7  
\put{$\bullet$} at 5.35 0.7 
\put{$\bullet$} at 5.45 0.7
\put{$\bullet$} at 5.55 0.7
\put{$\bullet$} at 5.65 0.7
\put{$\bullet$} at 5.75 0.7  
\put{$\bullet$} at 5.85 0.7
\put{$\bullet$} at 5.95 0.7
\put{$\bullet$} at 6.05 0.7
\put{$\bullet$} at 6.15 0.7
\put{$\bullet$} at 6.25 0.7  
\put{$\bullet$} at 6.35 0.7 
\put{$\bullet$} at 6.45 0.7
\put{$\bullet$} at 6.55 0.7
\put{$\bullet$} at 6.65 0.7
\put{$\bullet$} at 6.75 0.7  
\put{$\bullet$} at 6.85 0.7
\put{$\bullet$} at 6.95 0.7
\put{$\bullet$} at 7.05 0.7
\put{$\bullet$} at 7.15 0.7
\put{$\bullet$} at 7.25 0.7  
\put{$\bullet$} at 7.35 0.7 
\put{$\bullet$} at 7.45 0.7
\put{$\bullet$} at 7.55 0.7
\put{$\bullet$} at 7.65 0.7
\put{$\bullet$} at 7.75 0.7  
\put{$\bullet$} at 7.85 0.7
\put{$\bullet$} at 7.945 0.7 

\put{$R_2'$} at 0.3 0.37 
\put{$R_3'$} at 1.3 0.37 
\put{$R_4'$} at 1.3 2.07  
\put{$R_5'$} at 3.3 2.07  
\put{$R_6'$} at 3.3 4.37 
\put{$R_7'$} at 7.8 4.37  

\put{$x_0'$} at -0.08 -0.15 
\put{$x_1' \!\!=\!\! x_2'$} at 0.45 -0.15 
\put{$x_3' \!=\! x_4'$} at 1.45 -0.15 
\put{$x_5' \!=\! x_6'$} at 3.45 -0.15 
\put{$x_7' \!=\! x_8'$} at 7.95 -0.15 
\put{$y_0' \!=\! y_1'$} at -0.425 0 
\put{$y_2' \!=\! y_3'$} at -0.425 0.5 
\put{$y_4' \!=\! y_5'$} at -0.425 2.2 
\put{$y_6' \!=\! y_7'$} at -0.425 4.5 
\end{tiny} 

\put{} at -1.4 0 

\endpicture


\vspace{0.58cm}


We would like to close this section by giving a different type of choice 
for the sequence of rectangles which is simpler to describe and for which 
the preceding arguments are also valable. (For simplicity, we will use a 
similar construction to deal with the case $d > 2$, altough the preceding 
one still applies). This sequence $(R_m')_{m \geq 0}$ 
is of the form $[[0,x_m']] \times [[0,y_m']]$, where $(x_m')$ and 
$(y_m')$ are non decreasing sequences of positive integer numbers 
such that $x_0'=y_0'=0$, \esp $x_m' > x_{m-1}'$ and $y_m' = y_{m-1}'$ 
if $m$ is odd, and $x_m' = x_{m-1}'$ and $y_m' > y_{m-1}'$ if $m$ is even. 
If $(\ell_{i,j})$ is a double-indexed sequence of positive real numbers 
with total sum $\leq 1$, we chose these integer numbers in such a way 
that $x_{2m+1}' = x_{2m+2}' = [4^{m \tau_1}]$ and 
$y_{2m}' = y_{2m+1}' = [4^{m \tau_2}]$ for $m$ large enough. As before, 
inside the rectangle $R_m$ there is a ``good" vertical (resp. horizontal) 
segment of line $L_m$ for $m$ even (resp. odd). Therefore, for each 
$M_0 \in \mathbb{N}$ we can concatenate these segments between 
$L_{m-1} \cap L_m$ and $L_m \cap L_{m+1}$ at the 
$m^{\mathrm{th}}$ step for $m < M_0$, and between 
$L_{M_0 - 1} \cap L_{M_0}$ and the point of $L_{M_0}$ on the boundary 
of $R_{M_0}$ at the last step (see Figure 3). In this way we obtain a path 
(starting at the origin) of finite length $n(M_0)-1$ for which the sum 
$$\sum_{n = 0}^{n(M_0)-1} \ell_{i(n),j(n)}^{\tau_{\alpha(n)}}$$
is bounded by some number $S > 0$ which is independent of $M_0$.

Now let $f_k$, $k \!\in\! \{1,2\}$, be two commuting circle diffeomorphisms of class 
$C^{1 + \tau_k}$ which are not simultaneously conjugate to rotations. Fix again 
one of the maximal wandering open intervals for the dynamics, say $I$, and let 
$\ell_{i,j} = |f_1^if_2^j(I)|$. (Note that $\sum_{i,j} \ell_{i,j} \leq 1$.) 
The method above gives us a family of finite paths, and each of these 
paths determines uniquely a sequence of compositions. Remark however that 
there is a little difference here, since we allow the use of the inverses 
of $f_1$ and $f_2$. Therefore, in order to apply Lemma \ref{lema-cle}, we 
will need to consider now $\{f_1,f_1^{-1},f_2,f_2^{-1}\}$ as being our 
system of generators, and therefore we put $\tau = \max \{\tau_1,\tau_2\}$ 
and $C = \max \{ C_1,C_2,C_1',C_2' \}$, where $C_i$ (resp. $C_i'$) is a 
$\tau_i$-H\"older constant for the function $\log(f_i')$ (resp. $\log((f_i^{-1})')$). 
As in the previous proof, we need to verify that, for some $M_0 \in \mathbb{N}$, 
there exists a non trivial 
element in the sequence of compositions $(h_n)$ associated to 
its corresponding finite path which sends $I$ inside the $L$-neighborhood of 
itself, where $L := |I| / 2\exp(2^{\tau} C S)$. As before, for proving 
this it suffices to show that for every $N$ there exists $r \in \mathbb{N}$ 
such that one has $h_{r+i+1} = f_1 h_{r+i}$ for each $i \in \{0,\ldots,N-1\}$, 
or $h_{r+i+1} = f_2 h_{r+i}$ for each $i \in \{0,\ldots,N-1\}$. However, this 
last property is always satisfied if $M_0$ is big enough so that the number 
of points with integer coordinates in the line segment $L_{M_0}$ contained 
in $R_{M_0} \setminus R_{M_0 - 1}$ is greater than $N$. Note that it is in 
this last argument where we use the fact that we keep only finite sequences 
of compositions, altough our method combined with a diagonal type argument 
easily shows the existence of an infinite sequence for which the sum 
(\ref{acontrolar}) converges.


\subsection{The general case}

\hspace{0.5cm} In the case $d = 2$, the ``good" paths leading to the 
sequence of compositions which allows to apply Lemma \ref{lema-cle} were 
obtained by concatenating horizontal and vertical lines. When $d > 2$ we 
will need to concatenate lines in several (namely $d$) directions, and 
the geometrical difficulty for doing this is evident: in dimension bigger 
than 2, two lines in different directions do not necessarily intersect. 
To overcome this difficulty we will use the fact that, at each 
step ({\em i.e.} inside each rectangle), there is not only one 
finite path which is good, but this is the case for a ``large 
proportion" of finite paths. We first reformulate Lemma 2.1 
in this direction.

\vspace{0.1cm}

\begin{lem}{\em Let $\ell_{i,j}$ be positive real numbers, where $i \in \{1,\ldots,m\}$ 
and $j \in \{1,\ldots,n\}$. Assume that the total sum of the $\ell_{i.j}$'s is less 
than or equal to $1$. If $\tau$ belongs to $]0,1[$ and $A > 1$, then for a proportion 
of indexes $k \in \{1,\ldots,n\}$ greater than or equal to $(1 - 1/A)$ we have}
$$\sum_{i = 1}^{m} \ell_{i,k}^{\tau} \leq A \frac{m^{1-\tau}}{n^{\tau}}$$

\label{menos-basico}
\end{lem}

\noindent{\bf Proof.} As in the proof of Lemma 2.1, the mean value of the function 
\begin{equation}
k \mapsto \sum_{i=1}^{m} \ell_{i,k}^{\tau}
\label{funcion}
\end{equation}
is less than or equal to \esp $m^{1 - \tau} / n^{\tau}$. 
\esp The claim of the lemma then follows as a direct 
application of Chebychev's inequality: the proportion of points 
for which the value of (\ref{funcion}) is greater than this mean 
value times $A$ cannot exceed $1/A$. $\hfill\square$

\vspace{0.4cm}

Now let $(\ell_{i_1,\ldots,i_d})$ be a multi-indexed sequence of positive 
real numbers having total sum $\leq 1$, and let $\tau_1,\ldots,\tau_d$ be 
real numbers in $]0,1[$. Starting with $R_0 = [[0,0]]^d$, let us consider 
a sequence $(R_m)_{m \geq 0}$ of rectangles of the form \esp 
$R_m = [[0,x_{1,m}]] \times \cdots \times [[0,x_{d,m}]]$ \esp satisfying  
\esp $x_{k,m} \geq x_{k,m-1}$ \esp for each 
$k \in \{1,\ldots,d\}$, 
with strict inequality if and only if \esp $k \equiv m \esp (\mathrm{mod} \esp d)$. 
\esp For each $m \geq 1$ denote by $s(m) \in \{1,\ldots,d\}$ the residue class 
$(\mathrm{mod} \esp d)$ of $m$, and denote by $F_m$ the face  
$$[[0,x_{1,m}]] \times \cdots \times [[0,x_{s(m)-1,m}]] \times \{0\} \times 
[[0,x_{s(m)+1,m}]] \times \cdots \times [[0,x_{d,m}]]$$ 
of $R_m$. For each $(i_1,\ldots,i_{s(m)-1},0,i_{s(m)+1},\ldots,i_d)$ 
belonging to this face $F_m$ we consider the sum
$$\sum_{j = 0}^{x_{s(m),m}} 
\ell^{\tau_{s(m)}}_{i_1,\ldots,i_{s(m)-1},j,i_{s(m)+1},\ldots,i_d}.$$ 
By Lemma \ref{menos-basico}, if $A_m > 1$ then the proportion 
of points in $F_{m}$ for which this sum is bounded by 
$$A_m \cdot \frac{(1+x_{s(m),m})^{1-\tau_{s(m)}}}
{\prod\limits_{j \neq s(m)} (1+x_{j,m})^{\tau_{s(m)}}} 
= A_m \cdot \frac{X_{s(m),m}^{1-\tau_{s(m)}}}
{\prod\limits_{j \neq s(m)} X_{j,m}^{\tau_{s(m)}}}$$
is at least equal to \esp $(1-1/A_m)$, \esp 
where $X_{j,m} := 1 + x_{j,m}$. In order to concatenate 
the corresponding lines we will use the following elementary lemma.

\newpage


\beginpicture

\setcoordinatesystem units <1.13cm,1.13cm>

\putrule from 0 0 to 8 0 
\plot 8 0 7.8 0.1 /
\plot 8 0 7.8 -0.1 /
\put{$s(m)$-direction} at 9.2 0 

\putrule from 0 0 to 0 6 
\plot 0 6 -0.1 5.8 /
\plot 0 6 0.1 5.8 /

\putrule from 0 0.01 to 4 0.01 
\putrule from 0 -0.01 to 4 -0.01  

\putrule from 0.01 0 to 0.01 1.5  
\putrule from -0.01 0 to -0.01 1.5  

\putrule from 0 1.5 to 4 1.5 
\putrule from 4 0 to 4 1.5 
\putrule from 5 1 to 5 2.5 
\putrule from 7 3 to 7 4.5 
\putrule from 3 4.5 to 7 4.5 
\putrule from 1 2.5 to 5 2.5 

\plot 4 0 7 3 /
\plot 4.5 4.5 5.4 5.4 /
\plot 4 1.5 7 4.5 /
\plot 0 1.5 3 4.5 /

\plot 5.4 5.4 5.35 5.19 /
\plot 5.4 5.4 5.15 5.32 /

\put{$\bullet$} at 0 0.7 
\put{$C_m$} at -0.3 1.3 
\put{$C_{m+1}$} at 1 0.2 
\put{$R_m$} at 1.57 2.7 
\put{$R_{m+1}$} at 3.57 4.7  
\put{$F_m$} at 0.8 2  
\put{$F_{m+1}$} at 3.6 0.2 
\put{$\bullet$} at 2.3 0 
\put{$L_m$} at 4.37 1.45 
\put{$L_{m+1}$} at 4.4 3.4
\put{$s(m+1)$-direction} at 5.6 5.7 

\plot 2.285 0.7 5.285 3.7 /
\plot 2.3 0.7 5.3 3.7 /
\plot 2.305 0.704 5.302 3.704 /
\plot 2.315 0.7065 5.3 3.7065 /

\putrule from 0.55 1.25 to 4.55 1.25 
\putrule from 0.55 1.24 to 4.55 1.24
\putrule from 0.55 1.235 to 4.55 1.235

\putrule from -2.2 -0.5 to -2.2 0.35 
\putrule from -2.2 -0.5 to 0 -0.5 

\plot 0 -0.5 -0.2 -0.4 /
\plot 0 -0.5 -0.2 -0.6 /


\setdots
\plot 0 0 5.4 5.4 /
\putrule from 3 3 to 3 4.5 
\putrule from 3 3 to 7 3 
\putrule from 1 1 to 5 1 
\putrule from 1 1 to 1 2.5 

\putrule from 0 0.7 to 4 0.7 
\putrule from 1 1.7 to 5 1.7 

\plot 0 0.7 1 1.7 /
\plot 4 0.7 5 1.7 /

\plot 1 1.7 3 3.7 /
\plot 5 1.7 7 3.7 /
\putrule from 3 3.7 to 7 3.7 

\put{Figure 4} at 3 -1.2

\begin{tiny}
\put{$(i_1,\ldots,i_{s(m+1)-2},0,0,i_{s(m+1)+1},\ldots,i_d)$} at -2.33 0.85 
\put{admissible in $C_m$} at -2.2 0.55 

\put{$(i_1,\ldots,i_{s(m+1)-1},0,0,i_{s(m+1)+2},\ldots,i_d)$} at 2.3 -0.3 
\put{admissible in $C_{m+1}$} at 2.3 -0.6 
\end{tiny}

\put{} at -4.38 0 

\endpicture


\vspace{0.35cm} 

\begin{lem} {\em Let us chose inside each rectangle $(R_m)_{m \geq 1}$ a set 
$\mathcal{L}(m)$ of (complete) lines in the corresponding $s(m)$-direction whose 
proportion (with respect to all the lines in that direction inside $(R_m)$) is at 
least $(1-1/A_m)$. If $M_0 \!\in\! \mathbb{N}$ is such that 
$\sum_{m = 1}^{M_0} 1/A_m \!<\! 1$, then there exists a sequence of lines 
$L_m \in \mathcal{L}(m)$, $m \in \{0,\ldots,M_0\}$, such that $L_{m+1}$ 
intersects $L_m$ for every $m < M_0$.}
\label{clave-densidad}
\end{lem}

\noindent{\bf Proof.} Let us denote by $C_m$ the $(d-2)$-dimensional face of $R_m$ given by  
$$[[0,x_{1,m}]] \times \cdots \times [[0,x_{s(m)-1,m}]] \times \{0\} \times \{0\} 
\times [[0,x_{s(m)+2,m}]] \times \cdots \times [[0,x_{d,m}]].$$  
Call a point $(i_1,\ldots,i_{s(m)-1},0,0,i_{s(m)+2},\ldots,i_d) \in C_m$ {\em admissible} 
if there exists a sequence of lines 
$L_i \!\in\! \mathcal{L}(i), \esp i \!\in\! \{0,\ldots,m\}$, such 
that $L_i$ intersects $L_{i+1}$ for every $i \!\in\! \{0,\ldots,m-1\}$, and such that $L_m$ 
projects in the $s(m)$-direction into a point 
$(i_1,\ldots,i_{s(m)-1},0,i_{s(m)+1},i_{s(m)+2},\ldots,i_d) \in F_m$ for some 
$i_{s(m)+1} \!\in\! [[0,x_{s(m)+1,m+1}]]$. We will show that the proportion 
of admissible points in $C_{M_0}$ is greater than or equal to 
$$P := 1 - \sum_{m = 1}^{M_0} A_m > 0.$$

To prove this, for each $m \geq 0$ let us denote by $P_m$ the proportion of admissible points 
in $C_m$. Since $R_0$ reduces to the origin, it suffices to show that, for all $m \geq 0$, 
$$P_{m + 1} \geq P_m - \frac{1}{A_{m+1}}.$$
To prove this inequality first note that each line 
$L_{m+1} \in \mathcal{L}(m+1)$ determines 
uniquely a point $(i_1,\ldots,i_{s(m+1)-1},0,i_{s(m+1)+1},\ldots,i_d) \!\in\! F_{m+1}$. 
The projection into $C_m$ of this line then corresponds to the point 
$$(i_1,\ldots,i_{s(m+1)-2},0,0,i_{s(m+1)+1},\ldots,i_d).$$
If this is an admissible point of $C_m$ then we can concatenate the line $L_{m+1}$ 
to the sequence of lines corresponding to it (see Figure 4). Now the proportion 
of lines in $\mathcal{L}(m+1)$ being at least $1 - 1/ A_{m+1}$, the proportion of those lines 
which project on $C_m$ into an admissible point is at least equal to 
$$1 - \frac{1}{A_{m+1}} - (1 - P_m) = P_m - \frac{1}{A_{m+1}}.$$
By projecting in the $(s(m+1)+1)$-direction, this obviously implies that the 
proportion of admissible points in $C_{m + 1}$ is also greater than or equal 
to $P_m - 1/A_{m+1}$, thus finishing the proof. $\hfill\square$

\vspace{0.35cm}

Observe that a sequence of lines $L_m$ as above determines a finite 
path (starting at the origin) of points $(x_1(n),\ldots,x_d(n))$  
having non negative integer coordinates such that the distance 
between two consecutive ones is equal to $1$. Moreover, if we denote 
by $n(M_0)$ the length of this path plus $1$, the corresponding sum 
\begin{equation}
\sum_{n = 0}^{n(N_0)-1} \ell_{x_1(n),\ldots,x_d(n)}^{\tau_{\alpha(n)}}
\label{chito}
\end{equation}
is bounded by
\begin{equation}
\sum_{m = 0}^{M_0} A_m \cdot \frac{(1+x_{s(m),m})^{1 - \tau_{s(m)}}}
{\prod\limits_{i \neq s(m)} (1+x_{i,m})^{\tau_{s(m)}}} = 
\sum_{m = 0}^{M_0} A_m \cdot \frac{X_{s(m),m}^{1-\tau_{s(m)}}}
{\prod\limits_{j \neq s(m)} X_{j,m}^{\tau_{s(m)}}},
\label{cota}
\end{equation}
where $\alpha(n)$ equals the unique index in $\{1,\ldots,d\}$ for 
which $|x_{\alpha(n)}(n+1) - x_{\alpha(n)}(n)| = 1$. 

Now let us define $A_m \!=\! 2^{\varepsilon m \tau_{s(m)} /2} A$, where 
$A$ is a large enough constant so that $\sum_{m \geq 0} 1/A_m \!<\! 1$, 
and let us consider any choice of the $x_{k,m}$'s so that 
$X_{k,m}\simeq 2^{m\tau_k}$. For such a choice we have 
\begin{equation}
\frac{ X_{k,m}^{1-\tau_{k}}}{\prod\limits_{j \neq k}
X_{j,m}^{\tau_{k}}} = X_{k,m}^{-\varepsilon} \cdot \prod\limits_{j
\neq k} \frac{X_{k,m}^{\tau_j}}{X_{j,m}^{\tau_k}} \simeq
2^{-\varepsilon m \tau_k} \cdot \prod\limits_{j \neq k}
\frac{(2^{ m \tau_k })^{\tau_j}}{(2^{m \tau_j })^{\tau_k}} =
2^{-\varepsilon m \tau_k },
\label{chiito}
\end{equation}
where $\varepsilon := 1 - \tau_1 - \cdots - \tau_d > 0$. 
Therefore, for each 
$M_0 \in \mathbb{N}$ the preceding lemma provides us a 
sequence of lines $L_m$, $m \in \{0,\ldots,M_0\}$, such that 
$L_{m+1}$ intersects $L_m$ for each $m < M_0$, and such that 
the corresponding expression (\ref{cota}) is bounded from above by 
\begin{equation}
\label{chitito}
\sum_{m = 0}^{M_0} 2^{\varepsilon m \tau_{s(m) / 2}} A \cdot 
\frac{ X_{k,m}^{1-\tau_{k}}}{\prod\limits_{j \neq k} X_{j,m}^{\tau_{k}}}
\leq A C' \sum_{m \geq 0}  2^{-\varepsilon m \tau_{s(m)}/2 } \leq A
C' \sum_{m \geq 0}  2^{-\varepsilon m \tau'/2 } =:
S < \infty,
\end{equation}
where $\tau' := \min \{\tau_1,\ldots,\tau_d\}$ and $C'$ is a constant 
(independent of $M_0$) giving an upper bound for the quotient between 
the left and the right hand expressions in (\ref{chiito}). 

\vspace{0.4cm}

With all this information in mind we can proceed to the proof of Theorem A in the 
case $d > 2$ in the very same way as in the (second proof for the) case $d=2$. 
Indeed, assume that $f_k$, $k \in \{1,\ldots,d\}$, are circle diffeomorphisms 
as in the statement of the theorem which are not conjugate to rotations, and let 
$I$ be a maximal open wandering interval for the dynamics ({\em i.e.} a connected 
component of the complement of the minimal invariant Cantor set). 
Clearly, we can apply all our previous discussion to the multi-indexed sequence 
$(\ell_{i_1,\ldots,i_d})$ defined by  
$\ell_{i_1,\ldots,i_d} = |f_1^{i_1} \cdots f_d^{i_d} (I)|$. In particular, for 
each $M_0 \in \mathbb{N}$ we can find a finite path so that the sum (\ref{chito}) 
is bounded by the number $S > 0$ defined by (\ref{chitito}) 
(which is independent of $M_0$). Each such a path 
induces canonically a finite sequence of compositions by the $f_k$'s and their 
inverses. Therefore, in order to apply Lemma \ref{lema-cle} to get a contradiction, 
we need to verify that some of such sequences contains a (non trivial) element 
$h_n$ which sends $I$ into its $L$-neighborhood for $L:=|I|/2\exp(2^{\tau} C S)$, 
where $\tau:=\max \{ \tau_1,\ldots,\tau_d \}$ and 
$C:=\max \{C_1,\ldots,C_d,C_1',\ldots,C_d'\}$, with $C_k$ (resp. $C_k'$) being the  
$\tau_k$-H\"older constant of the function $\log(f_k')$ (resp. $\log((f_k^{-1})')$. 
To ensure this last property let $U$ be the $L$-neighborhood of $I$, and let 
$N \in \mathbb{N}$ be such that, given any wandering interval, among the first 
$N$ iterates of $f_1$, as well as for $f_2,\dots,f_d$, at least one of them 
sends this interval inside $U$. If we take $M_0$ large enough so that the 
number of points with integer coordinates in $L_{M_0}$ which are contained 
in $R_{M_0}\setminus R_{M_0-1}$ exceeds $N$, then one can easily see that 
the associated sequence of compositions contains the desired element 
$h_n$. This finishes the proof of Theorem~A.


\section{Proof of Theorem B}

\hspace{0.5cm} The strategy for the proof of Theorem B is well known. We prescribe 
the rotation numbers $\rho_1,\ldots,\rho_d$ (which are supposed to be independent 
over the rationals), we fix a point $p \in \clo$, and for each 
$(i_1,\ldots,i_d) \in \mathbb{Z}^d$ we replace the point 
$R_{\rho_1}^{i_1} \cdots R_{\rho_d}^{i_d}(p)$ 
by an interval $I_{i_1,\ldots,i_d}$ of length $\ell_{i_1,\ldots,i_d}$  
in such a way that the total sum of the $\ell_{i_1,\ldots,i_d}$'s is finite. 
Doing this we obtain a new circle on which the rotations $R_{\rho_k}$ 
induce nice homeomorphisms if we extend them apropiately to the 
intervals $I_{i_1,\ldots,i_d}$ (outside these intervals the 
induced homeomorphisms are canonically defined). More precisely, 
as it is well explained in \cite{DKN,herman,growth,TsP}, 
if there exists a constant $C' > 0$ so that for all 
$(i_1,\ldots,i_d) \in \mathbb{Z}^d$ and all $k \in \{1,\ldots,d\}$ one has  
\begin{equation}
\left| \frac{\ell_{i_1,\ldots,1+i_k,\ldots,i_d}}{\ell_{i_1,\ldots,i_k,\ldots,i_d}} - 1 \right| 
\frac{1}{\ell_{i_1,\ldots,i_k,\ldots,i_d}^{\tau_k}} \leq C',
\label{aver}
\end{equation}
then one can perform the extension to the intervals $I_{i_1,\ldots,i_d}$ 
in such a way the resulting maps $f_k$, $k \!\in\! \{1,\ldots,d\}$, are 
respectively $C^{1 + \tau_k}$ diffeomorphisms and commute, and moreover 
their derivatives are identically equal to $1$ on the invariant minimal 
Cantor set.\footnote{Condition (\ref{aver}) is also necessary under these 
requirements. Indeed, there must exist a point in $I_{i_1,\ldots,i_k,\ldots,i_d}$ 
for which the derivative of the corresponding map $f_k$ equals \esp 
$\ell_{i_1,\ldots,1+i_k,\ldots,i_d}/\ell_{i_1,\ldots,i_k,\ldots,i_d}$. 
\esp Since the derivative of $f_k$ at the end points of 
$I_{i_1,\ldots,i_k,\ldots,i_d}$ is assumed to be equal to $1$, 
condition (\ref{aver}) holds for $C^{\esp_{\! '}}$ being the $\tau_k$-H\"older 
constant of the derivative of $f_k$.} Indeed, one possible extension 
is given by $f_k (x) = (\varphi_{I_{i_1,\ldots,i_k,\ldots,i_d}})^{-1} 
\circ \varphi_{I_{i_1,\ldots,1+i_k,\ldots,i_d}} (x),$
where $x$ belongs to the interior of the interval $I_{i_1\ldots,i_k,\ldots,i_d}$. 
Here, $\varphi_I\!\!: ]a,b[ \rightarrow \mathbb{R}$ denotes the map 
$$\varphi_I (x) = \frac{-1}{b-a} \esp \mathrm{ctg} \Big( \pi \frac{x-a}{b-a} \Big).$$

It turns out that a good choice for the lengths is 
$$\ell_{i_1,\ldots,i_d} = \frac{1}{1 + |i_1|^{1/\tau_1} + \cdots |i_d|^{1/\tau_d}}.$$ 
Indeed, on the one hand, if we decompose the sum of the $\ell_{i_1,\ldots,i_d}$'s 
according to the biggest $|i_j|^{1 / \tau_j}$ we obtain
$$\sum_{(i_1,\ldots,i_d) \in \mathbb{Z}^d} \ell_{i_1,\ldots,i_d} \leq 
1 + \sum_{k = 1}^d \sum\limits_{\small\begin{array}{c} |i_j|^{1 / \tau_j} \leq |i_k|^{1 / \tau_k}\\
\mbox{ for all } j \in \{1,\ldots,d\}\\
|i_k| \geq 1 \end{array}} \frac{1}{1+|i_1|^{1/\tau_1}+\cdots |i_d|^{1/\tau_d}},$$
and therefore, for some constant $C > 0$, this sum is bounded by 
\begin{multline*}
1+\sum_{k=1}^d \sum_{n \geq 0} 
\frac{\mbox{card} \{ (i_1,\ldots,i_d) \!: |i_j|^{1 / \tau_j} \leq n^{1 / \tau_k} \mbox{ for all } 
j \!\in\! \{1,\ldots,d\}, i_k = n  \}}{1 + n^{1/\tau_k}} \\
\leq 1+C \sum_{k=1}^d \sum_{n \geq 1} \frac{1}{n^{1/\tau_k}} 
\prod\limits_{j \neq k} n^{\tau_j / \tau_k} 
= 1+C \sum_{k=1}^d \sum_{n \geq 1} 
\frac{ n^{(\sum_{j \neq k} \tau_j) / \tau_k}} {n^{1 / \tau_k}} \\
= 1+C \sum_{k=1}^{d} \sum_{n \geq 1} 
\frac{n^{(1-\tau_k - \varepsilon)/\tau_k}}{n^{1/\tau_k}} 
= 1+C \sum_{k=1}^{d} \sum_{n \geq 1} \frac{1}{n^{1 + \varepsilon / \tau_k}}, 
\end{multline*}
where $\varepsilon := 1 - (\tau_1 + \cdots + \tau_d)$. (Remark that, 
since $\varepsilon > 0$, the last infinite sum converges.)

On the other hand, the left hand expression in (\ref{aver}) is equal to 
\begin{multline*}
F(i_1,\ldots,i_d) := \left| \frac{|1+i_k|^{1/\tau_k} - |i_k|^{1/\tau_k}}
{1 + |i_1|^{1/\tau_1} + \cdots + |1+i_k|^{1/\tau_k} + \cdots + |i_d|^{1/\tau_d}} \right| \times \\
\times \Big( 1 + |i_1|^{1/\tau_1} + \cdots + |i_k|^{1/\tau_k} + \cdots + |i_d|^{1/\tau_d} \Big)^{\tau_k}.
\end{multline*}
In order to obtain an upper bound for this expression first note that, if $i_k \geq 0$, 
then 
$$F(i_1,\dots,i_k,\dots,i_d) \leq F(i_1,\dots,-1-i_k,\dots,i_d).$$
Therefore, we can restrict to the case where $i_k < 0$. For this case, 
denoting $B=1+\sum_{j\neq k} |i_j|^{1/\tau_j}$ and $a=|i_k|$ we have  
\begin{multline*}
F(i_1,\dots,i_d)=
\frac{a^{1/\tau_k}-(a-1)^{1/\tau_k}}{B+(a-1)^{1/\tau_k}}
\cdot \left(B+a^{1/\tau_k}\right)^{\tau_k}= \\
= \frac{a^{1/\tau_k}-(a-1)^{1/\tau_k}}{\left(B+(a-1)^{1/\tau_k}\right)^{1-\tau_k}}
\cdot \left(\frac{B+a^{1/\tau_k}}{B+(a-1)^{1/\tau_k}}\right)^{\tau_k}.
\end{multline*}
Both factors in the last expression are decreasing in $B$. 
Thus, since $B \geq 1$,
$$F(i_1,\ldots,i_d) \leq 
\frac{a^{1/\tau_k}-(a-1)^{1/\tau_k}}{\left(1+(a-1)^{1/\tau_k}\right)^{1-\tau_k}}
\cdot \left(\frac{1+a^{1/\tau_k}}{1+(a-1)^{1/\tau_k}}\right)^{\tau_k}.$$
Now note that $a \geq 1$. For $a=1$ the right hand expression above equals 
$2^{\tau_k}$. If $a>1$ then the Mean Value Theorem gives the estimate \esp 
\esp $a^{1/\tau_k} - (a-1)^{1/\tau_k} \leq a^{\frac{1}{\tau_k} - 1} / \tau_k$, 
\esp \esp and therefore the preceding expression is bounded from above by 
$$\frac{1}{\tau_k} \frac{a^{\frac{1}{\tau_k} - 1}}
{((a-1)^{1/\tau_k})^{1-\tau_k}} 
\cdot \left(\frac{a^{1/\tau_k}}{(a-1)^{1/\tau_k}}\right)^{\tau_k} = 
\frac{1}{\tau_k} \left(\frac{a}{a-1}\right)^{\frac{1}{\tau_k} -1}
\cdot \left(\frac{a}{a-1}\right) \le \frac{1}{\tau_k} \cdot
2^{\frac{1}{\tau_k} -1} \cdot 2 = \frac{2^{1/\tau_k}}{\tau_k}.$$
We have then shown that for any $(i_1,\dots,i_d)\in\mathbb{Z}^d$ one has
$$F(i_1,\dots,i_d) \le \frac{1}{\tau_k} 2^{1/\tau_k}.$$
In other words, if $\tau' = \min\{\tau_1,\ldots,\tau_d\}$ then inequality 
(\ref{aver}) holds for each $(i_1,\dots,i_d)\in\mathbb{Z}^d$ and every 
$k \in \{1,\ldots,d\}$ for the constant $C' = 2^{1/\tau'} / \tau'$, 
and this finishes the proof of Theorem B.

\vspace{0.15cm}


\begin{footnotesize}

\vspace{0.1cm}

Victor Kleptsyn\\

Universit\'e de Gen\`eve,
2-4 rue du Li\`evre, Case postale 64, 1211 Gen\`eve 4, Suisse 
(Victor.Kleptsyn@math.unige.ch)\\

\vspace{0.2cm}

Andr\'es Navas\\

Universidad de Santiago de Chile, 
Alameda 3363, Santiago, Chile (andnavas@uchile.cl)\\

\end{footnotesize}

\end{document}